\documentclass[12pt,a4paper,draft]{amsart}
\usepackage{epsf,amsfonts,amscd,latexsym}
\usepackage[T2A]{fontenc}
\usepackage[english]{babel}
\def\hepsffile{\leavevmode\epsffile}

\title[Fano variety of index 2 and degree 1]
{Non-rationality of a three-dimensional Fano variety of index 2 and degree 1.}
\thanks{This work was supported by the
grants RFBR no. 00--15--96085 and 02--01--00441,
INTAS no. 00---0259 and 00--0269, and NWO-RFBR no. 047-008-005}
\author{Mikhail Grinenko}
\address{Steklov Mathematical Institute}
\email{grin@mi.ras.ru}
\date{}

\newtheorem{theorem}{\sc Theorem}[section]
\newtheorem{proposition}[theorem]{\sc Proposition}
\newtheorem{lemma}[theorem]{\sc Lemma}

\makeatletter

\@addtoreset{equation}{section}
\newcommand{\l@abcd}[2]{\hbox to\textwidth{#1\dotfill #2}}

\makeatother

\newcommand*{\mybegintheorem}[1]{\begin{trivlist}\it%
  \item[\hspace{\labelsep}{\bf #1}]}
  \newcommand*{\myendtheorem}{\end{trivlist}}
  \newenvironment*{theorem*}{\mybegintheorem{Theorem.}}{\myendtheorem}
  \newenvironment*{proposition*}{\mybegintheorem{Proposition.}}{\myendtheorem}
  \newenvironment*{corollary*}{\mybegintheorem{Corollary.}}{\myendtheorem}
  \newenvironment*{definition*}{\mybegintheorem{Definition.}}{\myendtheorem}

\theoremstyle{remark}

\renewcommand{\phi}{\varphi}
\renewcommand{\epsilon}{\varepsilon}

\newcommand{\lra}{\longrightarrow}

\newcommand{\PT}{{{\mathbb P}^3}}
\newcommand{\PTw}{{{\mathbb P}^2}}
\newcommand{\POn}{{{\mathbb P}^1}}
\newcommand{\ZA}{{\mathbb Z}}
\newcommand{\QA}{{\mathbb Q}}
\newcommand{\RA}{{\mathbb R}}
\newcommand{\CA}{{\mathbb C}}
\newcommand{\PA}{{\mathbb P}}
\newcommand{\FA}{{\mathbb F}}
\newcommand{\mc}{\mathcal}
\newcommand{\mf}{\mathfrak}

\newcommand{\eqdef}{\stackrel{\rm def}{=}}
\newcommand{\ord}{\mathop{\rm ord}\nolimits}
\newcommand{\Supp}{\mathop{\rm Supp}\nolimits}
\newcommand{\mult}{\mathop{\rm mult}\nolimits}

\newcommand{\Pic}{\mathop{\rm Pic}\nolimits}
\newcommand{\Bir}{\mathop{\rm Bir}\nolimits}
\newcommand{\Aut}{\mathop{\rm Aut}\nolimits}

\newcommand{\Bas}{\mathop{\rm Bas}\nolimits}

\newcommand{\Center}{\mathop{\rm Center}\nolimits}

\newcommand{\Proj}{\mathop{{\bf Proj}\:}\nolimits}

\newcommand{\eps}{\varepsilon}

\begin{document}
\begin{abstract}
It is proved that a smooth three-dimensional Fano variety $U$ of index 2
and degree 1 (the double cone over the Veronese surface) has only the 
following Mori structures: $U$ itself, and a 
"two-dimensional family" of fibrations $U_l\to\POn$ in del Pezzo surfaces of 
degree 1 arising from blow-ups of curves $l$ of genus 1 and degree 1. In particular,
$\Bir(U)=\Aut(U)$, $U$ has no structures of conic bundle, and $U$ is not rational.
\end{abstract}

\maketitle

\section{Preliminaries.}
\label{sec1}

Let $U$ be a non-singular three-dimensional Fano variety of index 2 and degree 1.
It is often called "the double cone over the Veronese surface" by the reason given
in section \ref{sec3}. We suppose
$$
\Pic(U)=\ZA H,
$$
where $H$ is the class of ample divisors with conditions $\dim|H|=2$, $K_V=-2H$, 
and $H^3=1$. Let $C$ be an effective curve of (arithmetical) genus 1 and degree 1,
i.e., $p_a(C)=1$ and $C\circ H=1$, and $l$ the class of $C$ in $A^2(U)$.
Then $A^2(U)=\ZA l$, $[H]^2=l$, and there exists a 2-dimensional family 
${\mc P}$ of effective curves of this class. General elements of 
${\mc P}$ are non-singular, but ${\mc P}$ contains a 1-dimensional sub-family
of rational curves with nodes or cusps.

For any $C\in{\mc P}$, the linear system $|H-C|$ is a pencil consisting of
del Pezzo surfaces of degree 1. Since $C$ is the complete intersection of any
2 elements of $|H-C|$, the blow-up $\phi_C: U_C\to U$ has a natural projection
$\rho_C:U_C\to\POn$, which is a Mori fibration in del Pezzo surfaces of degree
1. Let us note that if $C$ has a node, the corresponding $U_C$ has the simplest 
double point (given by the equation $x^2+y^2+z^2+w^2=0$ in some local 
coordinates), and if $C$ has a cusp, $U_C$ has a unique singular point of type
$x^2+y^2+z^2+w^3=0$.

We recall that, by definition, a birational map $\chi:V\dasharrow V'$ between two 
Mori fibrations $\rho:V\to S$ and $\rho':V'\to S'$ defines the same Mori structure
(or that $V/S$ and $V'/S'$ are birational over base), if
\begin{itemize}
\item in case $\dim S=\dim S'>0$, there exists a birational map $\psi$ such 
       that $\rho'\circ\chi=\psi\circ\rho$,
\item in case $\dim S=\dim S'=0$ ($V$ and $V'$ are $\QA$-Fano varieties), 
     $V$ and $V'$ are isomorphic to each other (whence $\chi\in\Bir(V)$).
\end{itemize}

This determines an equivalence relation for Mori fibrations, and for any 
threefold $X$ we define its set of Mori (fiber) structures 
${\mathcal {MS}}(X)$ as
$$
{\mathcal {MS}}(X)=\left\{\mbox{Mori fibrations birational to $X$}\right\}
\big/ \left\{\mbox{birational over base}\right\}
$$
Remark that the set of Mori structures is a birational invariant.
We are ready to formulate the main result.

\begin{theorem}
\label{main_th}
Let $U$ be a non-sungular Fano threefold of index 2 and degree 1. Then
$$
{\mathcal {MS}}(U)=\left\{U\right\}\cup
\left\{\mbox{$U_C\to\POn$, where $C\in{\mc P}$}\right\}.
$$
In particular, $U$ has a unique model of Fano variety and has no structures
of conic bundles, $\Bir(U)=\Aut(U)$, and $U$ is not rational. In general, 
$\Bir(U)\cong \ZA_2$.
\end{theorem}

Some attempts to study birational geometry of $U$ were undertaken in
\cite{Isk}, \cite{Kh}, and \cite{Grin1}, using the maximal singularity method.
In \cite{Isk} it is shown that linear systems on $U$ may have maximal 
singularities only at curves from ${\mc P}$, and some cases of
maximal singularities over points are dealt. Note that $|H-C|$ for 
$C\in{\mc P}$ gives an example of a linear system without fixed components
with a maximal singularity along $C$. The author of \cite{Kh}
tried to complete the work but failed in excluding maximal singurarities
over singular points of curves from ${\mc P}$ (see \cite{Grin1} for 
details). In \cite{Grin1}, the previous results were re-proved in a more regular
way, and finally the following assertion was achieved: if a linear system ${\mc D}$
without fixed components has a maximal singularity at $B_0$, then either
$B_0$ is a curve from ${\mc P}$, or $B_0$ is the singular point of
some $C\in{\mc P}$ for the corresponding maximal valuation on $U$
given by a weighted blow-up of type $(1,1,N)$, $2\le N\le 5$, in nodal 
case and  $(1,1,2)$ in cuspidal case (theorem 3.1 in \cite{Grin1}).

In order to complete proving theorem \ref{main_th}, we have to exclude 
maximal singularities over points of $U$ and, if ${\mc D}$ has a maximal 
singularity along $C\in{\mc P}$, jump onto $U_C/\POn$. On $U_C$, one needs to show
that the strict transform of ${\mc D}$ must not have super-maximal singularities,
and then we are done. 

The paper is organized as follows. In section \ref{sec2}, we recall some recent
results of M.Kawakita that allow to simplify much the usual argumentation in the
maximal singularities method. Section \ref{sec3} contains the description of geometry
of $U$ and $U_C/\POn$. Section \ref{sec4} deals with cases $(1,1,N)$ over points
mentioned above. We exclude super-maximal singularities over non-sungular points 
of $U_C/\POn$ in section \ref{sec5}, and over singular points in section 
\ref{sec6}. The last section contains the conclusion.

We assume that all required statements of the maximal singularities method
(see \cite{Pukh} for reference) and the Sarkisov program (\cite{Corti})
are known.

\noindent{\bf Acknowledgment.} The present paper was written during my stay
at the Max-Planck Institute for Mathematics in Bonn. I am very grateful to
the Institute for the kind support and hospitality.

\section{Maximal singularities and weighted blow-ups.}
\label{sec2}

Let ${\mc D}$ be a linear system without fixed components on a threefold $V$.
We assume $V$ to be in the Mori category. Suppose that ${\mc D}$ has a maximal 
singularity, i.e., there exists a discrete valuation ${\mf v}$ centered
at a point or a curve on $V$ such that the log pair $K_V+\frac1m{\mc D}$ is not
canonical with respect to ${\mf v}$, where $m$ is the adjunction threshold 
of ${\mc D}$. The argumentation of \cite{Corti} shows that there exists a 
divisorial contraction $\phi:V'\to V$ in the Mori category such that the 
exceptional divisor $E'$ of $\phi$ also defines a maximal singularity of 
${\mc D}$. The corresponding discrete valuation ${\mf v}_{div}$ has two
main properties. First, $\Center_V{\mf v}_{div}\subset \Center_V{\mf v}$.
Second, we may suppose that the minimal discrepancy of the log pair 
$K_V+\frac1m{\mc D}$ corresponds to ${\mf v}_{div}$. In other words, if 
$\psi:W\to V$ is a log resolution with exceptional divisors 
$\{E_1,E_2,\ldots\}$, and $E_1$ corresponds to ${\mf v}_{div}$, then
in the equation
$$
K_W+\frac1m\psi_*^{-1}{\mc D}=\psi^*\left(K_V+\frac1m{\mc D}\right)+
\sum a\left(W,K_V+\frac1m{\mc D},E_i\right)E_i
$$
we have $a(W,K_V+\frac1m{\mc D},E_1)\le a(W,K_V+\frac1m{\mc D},E_i)$ for all $i$.
Brief explanations can be also found in section 2 of \cite{Grin1}.

Recently the following remarkable results were achieved (\cite{Kaw1}, 
theorem 2.2, and \cite{Kaw2}, theorem 2.5):
\begin{theorem}
\label{Kaw_th}
Let $X$ be the germ of a three-dimensional variety, and $\phi:Y\to X$ is a 
divisorial contraction in the Mori category to a point $B_0\in X$. Then:

\noindent i) if $B_0$ is a smooth point of $X$, $\phi$ is the weighted blow-up
with weights $(1,L,N)$ in suitable local coordinates $[x,y,z]$ on $X$, where
$L\le N$ are coprime integers;

\noindent ii) if $B_0$ is a singular point of $X\in\CA^4$ given by 
$xy+z^2+w^2=0$ in suitable local coordinates, $\phi$ is the weighted blow-up
with $wt(x,y,z,w)=(1,1,1,1)$;

\noindent iii) if $B_0$ is a singular point of $X\in\CA^4$ given by 
$xy+z^3+w^2=0$ in suitable local coordinates, $\phi$ is the weighted blow-up
with $wt(x,y,z,w)=(1,1,1,1)$ or $wt(x,y,z,w)=(1,5,2,3)$.
\end{theorem}

In other words, if ${\mc D}$ has a maximal singularity over a point 
$B_0\in V$, we can always assume that there exists a maximal singularity
over $B_0$ which is realized as a divisorial contraction in the Mori
category, and if $B_0$ is one of the listed above, this contraction
is inverse to the corresponding weighted blow-up.

The assertion {\it i) } of theorem \ref{Kaw_th} has the following 
interpretation. Let ${\mf v}$ be any geometrical discrete 
valuation centered at a smooth point $B_0\in V$, and 
$\phi_1:V_1\to V$ be the blow-up of $B_0$ with the exceptional 
divisor $E_1\subset V_1$, $E_1\cong\PTw$. Denote by $B_1\subset E_1$ the 
center of ${\mf v}$ on $V_1$. If $B_1$ is also a point, we 
blow up $V_1$ at $B_1$ with the exceptional divisor $E_2\subset V_2$, and 
so on, until for some $L\ge 1$ we obtain $B_L$ is a curve on 
$E_L\subset V_L$. In the maximal singularities method we can always 
suppose that $B_L$ is a line on $E_L\cong\PTw$. Then blowing up $B_L$, we
obtain $\phi_{L+1}:V_{L+1}\to V_L$ with the exceptional divisor
$E_{L+1}\cong\FA_2$. Now $B_{L+1}$ has to be a section of
$E_{L+1}$ that covers $B_L$, and we proceed with this way, each time 
blowing up sections of exceptional divisors that are ruled surfaces. 
Finally, we must stop at $E_N$ which realizes ${\mf v}$.

We introduce the oriented graph $\Gamma({\mf v})$ as follows. 
It consists of $N$ vertices $\{1,2,\ldots,N\}$, and there is an arrow
$j\to i$, if $j>i$ and $B_{j-1}\subset E_i^{j-1}$, where leading 
indices mean the strict transform on the corresponding $V_i$'s. Note that
there are always arrows $i+1\to i$ since $B_i\subset E_i$ by choosing
the resolution of ${\mf v}$.

{\it The statement {\it i)} of theorem \ref{Kaw_th} means that 
$\Gamma({\mf v}_{div})$ is a chain, i.e., only arrows $i+1\to i$ 
exist. Moreover, the numbers $L$ and $N$ are the same as in the theorem.
Finally, for $i>L$, $B_i$ does not intersect the minimal section of
the ruled surface $E_i$. }

From the last assertion we can deduice the following proposition:

\begin{lemma}
\label{Ei_lem}
For $L<i<N$, we have $E_i\cong\FA_{i+1-L}$, 
${\mc N}_{E_i|V_i}={\mc O}_{E_i}(-s_i-f_i)$, $E_{i-1}^i|_{E_i}=s_i$,
and ${\mc N}_{B_i|V_i}\cong{\mc O}(-1)\oplus{\mc O}(i+1-L)$,
where $s_i$ and $f_i$ are the classes of the minimal section and a fiber
of $E_i$. Viewing $B_i$ as a curve on $E_i$, we obtain 
$B_{i}\sim s_i+(i+1-L)f_i$.
\end{lemma}
\noindent{\bf Proof.} Direct calculation, starting from $E_{L+1}$. Use the
fact that $B_i$ does not intersect the minimal section of $E_i$. \qed

The valuation ${\mf v}_{div}$ defines the numbers
$\nu_i=\mult_{B_{i-1}}{\mc D}^{i-1}$, $i=1,\ldots,N$. Since 
$K_V+\frac1m{\mc D}$ is not canonical with respect to ${\mf v}$, we
obtain the N\"other-Fano inequality in the form
$$
\nu_1+\nu_2+\ldots+\nu_N>m(L+N).
$$
Recall that $m$ is the adjunction threshold. For example, if 
${\mc D}\subset|nH|$ is a linear system on the double cone over the Veronses 
surface $U$, then $m=\frac{n}2$, and we have 
$\nu_1+\ldots+\nu_N>\frac{n}2(L+N)$.

\section{Geometric constructions.}
\label{sec3}

\noindent{\bf Equation and construction of $U$.} It is easy to observe that
$U$ can be realized as a hypersurface of degree 6 in the weighted ptojective
space $\PA=\PA(1,1,1,2,3)$. Let $[p,x,y,z,w]$ be the coordinates in $\PA$ 
with $wt(p,x,y,z,w)=(1,1,1,2,3)$. Since $U$ does not pass through the singular
points $(0,0,0,1,0)$ and $(0,0,0,0,1)$ of $\PA$, the equation of $U$ has
to contain monomials $w^2$ and $z^3$. Suppose that $B_0=(1,0,0,0,0)$, and $B_0$
is the singular point of $l_0\in{\mc P}$. We assume that the "hyperplanes"
$T=\{x=0\}$ and $\{y=0\}$ cut out $l_0$, and $T$ is tangent to $U$ at $B_0$.
To avoid complicating the notation, we keep $[x,y,z,w]$ as the corresponding 
coordinates in the affine part $\{p\ne 0\}$ of $\PA$.

Suppose that $B_0$ is the node of $l_0$. Then the equation of $U$ can be chosen in
the form
$$
w^2+z^3+z^2f_2(x,y)+zf_4(x,y)+x+f_6(x,y)=0,
$$
where $\deg f_i\le i$, $f_2(0,0)\ne 0$, $f_4(0,0)=f_6(0,0)=0$, and the minimal
degree of monomials in $f_6$ is not less than 2 (i.e., $f_6$ does not contain 
a linear part).

If $B_0$ is the cusp on $l_0$, we may assume the equation of $U$ to be in the form
$$
w^2+z^3+zf_4(x,y)+x+f_6(x,y)=0,
$$
where $f_4$ and $f_6$ satisfy the same requirements as before.

Being restricted to $U$, the projection $\PA\dasharrow\PA(1,1,1,2)$ gives
a morphism of degree 2 $U\stackrel{2:1}{\lra}\PA(1,1,1,2)$ branched over
a surface of degree 6. We can view $\PA(1,1,1,2)$ as the cone in
$\PA^6$ over the Veronese surface in $\PA^5$. This is why $U$ is often called
"the double cone over the Veronese surface".

\noindent{\bf Construction of the corresponding del Pezzo fibrations.} Let 
$C\in{\mc P}$, and $V=U_C\to U$ the blow-up of $C$. Recall that $C$ is the 
complete intersection of any two elements from $|H-C|$. The natural 
projection $\rho:V\to\POn$ realizes $V$ as a Mori fibrations on del Pezzo
surfaces of degree 1. We see that $V$ is always Gorenstein, all fibers of
$\rho:V\to\POn$ are irreducible and reduced, and
$$
\Pic(V)=\ZA[-K_V]\oplus\ZA F,
$$
where $F$ is the class of a fiber. $V$ is non-singular if $C$ is non-singular,
$V$ has an ordinary double point in the fiber corresponding to $T\in|H-C|$ 
(as before, $T$ is the tangent section of $U$ at $B_0$)
if $C$ has a node, and $V$ has a double point with a local equation
$ab+c^2+d^3=0$ if $C$ has a cusp.

Another construction of $\rho:V\to\POn$ can be achieved as follows (see
\cite{Grin2}, section 2.1, the case $n_1=0$, $n_2=1$, $n_3=2$.).
Let $X=\Proj_{\POn}{\mc O}\oplus{\mc O}\oplus{\mc O}(1)\oplus{\mc O}(2)$
with the natural projection $\pi:X\to\POn$. Denote by $M$ the class of the 
tautological bundle on $X$, by $L$ the class of a fiber, by $t_0$ the class of
minimally twisted sections (i.e., $t_0$ has an effective representative and
$t_0\circ M=0$), and by $p$ the class of a line in a fiber of $\pi:X\to\POn$.
Consider a hypersurface $Q\sim 2M-2L$ which is fibered over $\POn$ into 
non-degenerated two-dimensional quadric cones. Denote by $t_b$ the section
of $X/\POn$ that consists of the vertices of the cones. The class of $t_b$
is exactly $t_0$. Let $R_Q$ be the restriction of some effective divisor
$R\sim 3M$ to $Q$. Finally, take the double cover 
$\phi:V\stackrel{2:1}{\lra}Q$ branched over $R_Q$. We obtain a variety fibered
into del Pezzo surfaces of degree 1 over $\POn$. It is not very difficult to
see that $V\to\POn$ is isomorphic to $U_C$ for some $C\in{\mc P}$ under 
the suitable choice of $R_Q$, and $\rho=\pi\circ\phi$.

Let $s_b=\phi^{-1}(t_b)$. Then $s_b$ is a section of $\rho$, and 
$s_b$ cuts out the base point of the anticanonical linear system on each 
fiber. It can be defined also as
$$
s_b=\Bas |-K_V+lF|
$$
for $l\gg 0$. We define $s_0=\frac12\phi^*(t_0)$ and $f=\frac12\phi^*(p)$.
It is easy to see that $s_b\sim s_0$. Note that 
$$
NE(V)=\overline{NE(V)}=\RA_+s_0\oplus\RA_+f.
$$ 
Denote by $G_V$ the unique effective element of $|-K_V-2F|$. The divisor
$G_V$ is exceptional for the blow-up $V=U_C\to U$ and isomorphic to
$C\times\POn$. In \cite{Grin2} it was shown that $K_V^2\sim 2s_0+3f$ and
$s_0\circ(-K_V)=1$. Clearly, $f\circ F=0$ and $s_0\circ F=f\circ(-K_V)=1$.

\noindent{\bf Relations between divisors on $U$ and $V=U_C\to\POn$.}
Let $\phi_C:V=U_C\to U$ be the blow-up of $C\in{\mc P}$. Suppose that
${\mc D}\subset|aH-\mu C|$. In other words, ${\mc D}$ consists of divisors
of the class $aH$ that have the multiplicity $\mu$ along $C$. Let 
${\mc D}_V=(\phi_C)_*^{-1}\subset|n(-K_V)+mF|$ be the strict transform of
${\mc D}_U$ on $V$. Then it is easy to compute that
$$
\begin{array}{l}
n=a-\mu, \\
m=2\mu-a,
\end{array}
$$ 
i.e., ${\mc D}_V\subset|(a-\mu)(-K_V)+(2\mu-a)F|$.

\section{Excluding infinitely near singularities on $U$.}
\label{sec4}

Let $l_0\in{\mc P}$, $B_0$ the singular point of $l_0$. As it follows from 
theorem 3.1.2 of \cite{Grin1}, we have to consider the case of infinitely
near maximal singularities over $B_0$ that are realized by the weighted 
blow-up with weights either $(1,1,N)$, $2\le N\le 5$, if $B_0$ is a node, or 
$(1,1,2)$, if $B_0$ is a cusp. In all these cases we assume that 
$\mu=\mult_{l_0}{\mc D}\le\frac{n}2$, where ${\mc D}\subset|nH|$ is a linear
system without fixed components.

Following to section \ref{sec2}, weighted blow-ups of type $(1,1,N)$ can be 
resolved as follows. We blow up the point $B_0$, and the exceptional divisor
(say, $E_1$) contains the next center $B_1$, which must be a line. We blow up
$B_1$ with a ruled surface $E_2$ as the exceptional divisor, and this is all if 
$N=2$. Otherwise, there exists a section $B_2\subset E_2\setminus E_1^2$ which is
the next center. We blow up it, and so on. As in the end of section \ref{sec2},
we denote $\nu_i=\mult_{B_{i-1}}{\mc D}^{i-1}$, $\nu_1\ge\nu_2\ge\ldots>\frac{n}2$. 
The corresponding N\"other-Fano inequalities are given by
\begin{equation}
\label{NF1_ineq}
\nu_1+\ldots+\nu_N> \frac{n}2\left(1+N\right).
\end{equation}

\noindent{\bf The nodal case.} As it is shown in section \ref{sec3}, we may
choose local coordinates $[x,y,z,w]$ such that the local equation of $U$
near $B_0=(0,0,0,0)$ is
$$
w^2+z^3+z^2f_2(x,y)+zf_4(x,y)+x+f_6(x,y)=0,
$$
where $\deg f_i\le i$, $f_2(0,0)\ne 0$, $f_4(0,0)=f_6(0,0)=0$, $f_6$ does
not contain a linear part, and $T=\{x=0\}$ is the tangent space to $U$ at 
$B_0$. We set $T_U=T|_U\in|H-B_0|=|H-l_0|$.

\begin{figure}[htbp]
\begin{center}
\epsfxsize 10cm
\hepsffile{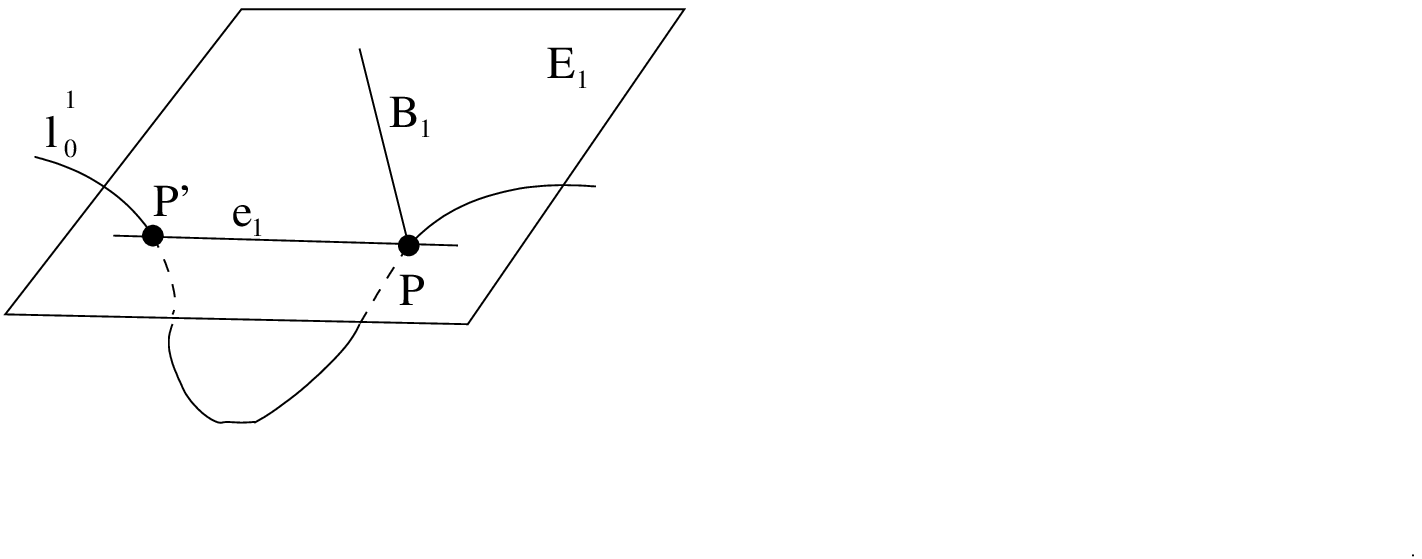}
\caption{}\label{ver1.fig}
\end{center}
\end{figure}

Let $\phi_1:U_1\to U$ be the blow-up of $B_0$, $E_1\cong\PTw$ the 
exceptional divisor, $l_0^1$ the strict transform of $l_0$ on $U_1$ 
(leading indices usually mean the strict transform on the corresponding
floor of the chain of blow-ups), $\{P,P'\}=l_0^1\cap E_1$ (clearly, $P\ne P'$),
$e_1\subset E_1$ the line joining $P$ and $P'$. The situation is given by
figure \ref{ver1.fig}. Note that except for $T_U$, the strict transform of any 
element $S\in|H-l_0|$ on $U_1$ cuts out the line $e_1$ on $E_1$. 

It is easy to observe from the equation of $U$ that the singular point of 
$T_U$ at $B_0$ is exactly $A_k$, $1\le k\le 5$. Respectively, the restriction 
$T_U^1|_{E_1}$ is either a non-singular conic that passes through $P$ and $P'$, 
or a couple of lines that contain $P$ and $P'$ and have their common point
outside of $e_1$ (thus they are different from $e_1$). In the latter case 
$T_U^1$ may have a singular point in the common point of the couple of  
the lines, but anyway $T_U^1$ is non-singular along $l_0^1$.

The possible situation is exhausted by the following cases:
\begin{itemize}
\item[\it A).] neither $P$ nor $P'$ is contained in $B_1$;
\item[\it B).] $e_1=B_1$;
\item[\it C).] $P\in B_1\ne e_1$ and $B_1\subset\Supp T_U^1|_{E_1}$;
\item[\it D).] $P\in B_1\ne e_1$, $T_U^1|_{E_1}$ is a conic, and $B_1$ is 
                 tangent to $T_U^1|_{E_1}$ at $P$;
\item[\it E).] other cases when $P\in B_1\ne e_1$.
\end{itemize}
  
Note that in all cases except for $e_1=B_1$ we have the inequality
\begin{equation}
\label{mu_ineq}
\nu_1\ge\nu_2+\mu,
\end{equation}
which arises from the intersection of a general line in $E_1$ through the point
$P'$ and an element from ${\mc D}$. Now we deal separately with the listed cases.

\medskip

\noindent{\it The case C).} Consider the family of divisors 
${\mc L}=\{z=ay, a\ne 0\}$. Clearly, ${\mc L}\subset|2H|$ with non-singular $K3$
surfaces as general elements. Then,  ${\mc L}^1$ has no base curves on $U_1$, 
so its elements cover $E_1$. Choose general elements $S\in{\mc L}$ and 
$D\in{\mc D}$, and consider the curve $C=T_U\cap S$. Clearly, $C\sim 2l$,
$C\not\subset D$, $C$ has a double point at $B_0$, and $C^1\cap B_1\ne\emptyset$.
Thus
$$
2n=D\circ C\ge 2\nu_1+\nu_2.
$$
Since $\nu_1+\ldots+\nu_N>n\left(1+\frac{N}2\right)$ from the N\"other-Fano
inequality, we see that $N=2$ or $N=3$ are impossible. In order to exclude
the cases $N>3$, apply the inequality \ref{mu_ineq}. We have
$$
2n\ge 2\nu_1+\nu_2\ge 2\mu+3\nu_2>2\mu+\frac32n,
$$
hence $\mu<\frac{n}4$.

Now let us proceed with the 
way which we will use repeatedly in the sequel. Choose general elements
$S\in|H-l_0|$ and $D\in{\mc D}$. $S$ is a non-singular del Pezzo surface of
degree 1. Set
$$
D|_S=\mu l_0+C,
$$
where the components of the residual curve $C$ do not contain $l_0$. We resolve 
the maximal valuation over $B_0$ in the way given in section \ref{sec2}.
Thus, we have the blow-ups $\phi_i:U_i\to U_{i-1}$, $U_0=U$, of the centers
$B_{i-1}\subset E_{i-1}$. In our case, $B_0$ is a point, the other $B_i$'s are
curves. As it has been mentioned, the corresponding oriented graph has no
incidencies and $B_i\cap E_{i-1}^i=\emptyset$ for $i\ge 2$. We put 
$e_i=E_i|_{S^i}$, $\tilde B_i=B_i\cap S^i$ for $i\ge 1$, and
$\tilde\nu_i=\mult_{\tilde B_{i-1}}D^{i-1}|_{S^{i-1}}$. We have
$$
\begin{array}{l}
D|_S=\mu l_0+C, \\
D^1|_{S^1}=\mu l_0^1+C^1+m_1e_1, \\
\cdots \\
D^{N-1}|_{S^{N-1}}=\mu l_0^{N-1}+m_1e_1^{N-1}+\ldots+m_{N-1} e_{N-1},
\end{array}
$$
where the numbers $m_i$ satisfy the relations 
\begin{equation}
\label{nu_m_eq}
\tilde\nu_i=nu_i+m_i.
\end{equation}
Note that $\tilde B_i\not\in e_{i-1}^i$. Set 
$\alpha_i=\mult_{\tilde B_i}C^i$, $\alpha_0\ge\alpha_1\ge\ldots$, and  let
$$
k=\max\{i: \tilde B_{i-1}\in l_0^{i-1}\}.
$$
Then we have
$$
\begin{array}{l}
\tilde\nu_1=2\mu+\alpha_0, \\
\tilde\nu_2=\mu+\alpha_1+m_1, \\
\cdots \\
\tilde\nu_k=\mu+\alpha_{k-1}+m_{k-1}.
\end{array}
$$
Moreover, one deduices that
$$
2\alpha_0+\alpha_1+\ldots+\alpha_{k-1}\le (D|_S-\mu l_0)\circ l_0=n-\mu.
$$
Suppose that $k<N$. Then $\tilde\nu_{k+1}=\alpha_k+m_k$. Thus
$$
\tilde\nu_1+\ldots+\tilde\nu_{k+1}=(k+1)\mu+\sum_{i=0}^k\alpha_i+
   \sum_{i=1}^km_i.
$$
Using (\ref{nu_m_eq}) and the inequality
$$
\sum_{i=0}^k\alpha_i\le 2\alpha_0+\alpha_1+\ldots+\alpha_{k-1}\le n-\mu,
$$
we obtain
$$
\nu_1+\ldots+\nu_{k+1}\le (k+1)\mu+n-\mu=n+k\mu.
$$
Since $\mu\le\frac{n}2$, we get 
$$
\nu_1+\ldots+\nu_{k+1}\le \frac{n}2(k+1),
$$
which is impossible because $\nu_i>\frac{n}2$ for all $i$. Hence $k=N$. Arguing as 
before, we see that
$$
\nu_1+\ldots+\nu_N\le (N+1)\mu+n-\mu=n+N\mu.
$$
Using the estimation $\mu<\frac{n}4$, one has
$$
\nu_1+\ldots+\nu_N< n\left(1+\frac{N}4\right).
$$
It only remains to compare it with the N\"oether-Fano inequality 
(\ref{NF1_ineq}), and we get the contradiction.

\medskip

\noindent{\it The case B).} In this case $e_1=B_1$. We repeat the argumentation
of the previous case with minor changes. Let ${\mc L}=|H-l_0|$,
$\phi_1:U_1\to U$ the blow-up of $B_0$. For general $S\in{\mc L}$, 
$S^1|_{E_1}=e_1$. We set
$$
\begin{array}{l}
D|_S=\mu l_0+C, \\
D^1|_{S^1}=\mu l_0^1+C^1+m_1 e_1,
\end{array}
$$
where $m_1=\nu_2$ since $e_1=B_1$. Let $\phi_2:U_2\to U_1$ be the blow-up of $B_1$,
$E_2\cong\FA_2$ the exceptional divisor with $s_2$ and $f_2$ as the classes
of the minimal section and a fiber respectively. Note that $E_1^2|_{E_2}$ cuts
out the minimal section of $E_2$, the pencil ${\mc L}^2|_{E_2}\subset|s_2+2f_2|$ 
has two base points, say, $Q$ and $Q'$, outside of the minimal section. 

We resolve the valuation as in section \ref{sec2}. Let
$$
k=\max\left\{i: B_{i-1}\cap l_0^{i-1}\ne\emptyset\right\}.
$$
Clearly, $k\ge 3$. For $i\ge 2$, let $\tilde B_i$ be any of the intersection
points $S^i\cap B_i$. We see that $\tilde B_k\not\in l_0^k$. Then, for $i\ge 2$,
$$
D^i|_{S^i}=\mu l_0^i+C^i+m_2 e_2^i+m_3 e_3^i+\ldots+m_i e_i,
$$
where $e_i=E_i|_{S^i}$. The important observation is that $m_2=0$.

Suppose that $k<N$. We have
$$
\begin{array}{l}
\nu_1+\nu_2=2\mu+\alpha_0, \\
\tilde\nu_3=\mu+\alpha_2, \\
\tilde\nu_4=\mu+\alpha_3+m_3,\\
\cdots \\
\tilde\nu_k=\mu+\alpha_{k-1}+m_{k-1},\\
\tilde\nu_{k+1}=\alpha_k+m_k,
\end{array}
$$
where $\tilde\nu_i=\mult_{\tilde B_{i-1}}D^{i-1}|_{S^{i-1}}=\nu_i+m_i$ and
$\alpha_i=\mult_{\tilde B_i}C^i$. Since $\mu\le\frac{n}2$ and
$\alpha_0+\alpha_2+\ldots+\alpha_k\le n-\mu$, we get a contradiction:
$$
\frac{n}2(k+1)<\nu_1+\ldots+\nu_{k+1}\le n+(k-1)\mu\le\frac{n}2(k+1).
$$ 
Suppose that $k=N$. Arguing as before, we see that
$$
\nu_1+\ldots+\nu_N\le N\mu+\alpha_0+\alpha_2+\ldots+\alpha_{N-1}\le
n+(N-1)\mu\le\frac{n}2(N+1),
$$
and this contradicts to the N\"other-Fano inequality (\ref{NF1_ineq}).

\medskip

\noindent{\it The case A).} As it has been mentioned above, in this case 
we have $\nu_1\ge\nu_2+\mu$. Denote $\tilde B_1=B_1\cap e_1$, 
$\tilde B_1$ is different from both $P$ and $P'$. We set
$$
\begin{array}{l}
D|_S=\mu l_0+C, \\
D^1|_{S^1}=\mu l_0^1+C^1+m_1 e_1.
\end{array}
$$
Since $\tilde B_1\not\in l_0^1$, we have
$$
\begin{array}{l}
\tilde\nu_1=2\mu+\alpha_0, \\
\tilde\nu_2=\alpha_1+m_1,
\end{array}
$$
and $\alpha_0+\alpha_1\le n-\mu$ (we keep the notation of the previous cases).
Thus,
$$
2\nu_2+\mu\le\nu_1+\nu_2\le 2\mu+\alpha_0+\alpha_1\le n+\mu,
$$ 
whence $n<2\nu_2\le n$, a contradiction.

\medskip

\noindent{\it The case E).} It is easy to compute that
$$
{\mc N}_{l_0^1|U_1}\cong{\mc O}(-1)\oplus{\mc O}(-3).
$$
Let $\phi_2:U_2\to U_1$ be the blow-up of $l_0^1$ (not $B_1$!!), $E_2\cong\FA_2$ 
the exceptional divisor with $s_2$ and $f_2$ as the classes of the minimal section
and a fiber. We see that
$$
D^2|_{E_2}\sim \mu s_2+(3\mu+n-2\nu_1)f_2.
$$
By the condition, $T_U^1|_{E_1}$ and $B_1$ are transversal to each other at $P$.
Since $T_U^2$ cuts out the minimal section of $E_2$, we deduice that $B_1^2$
intersects $E_2$ outside of the minimal section, and at this point $D^2|_{E_2}$
has the multiplicity not less than $\nu_2$. Thus
$$
3\mu+n-2\nu_1\ge\nu_2,
$$
and since $\nu_1\ge\nu_2+\mu$, we get
$$
n+\mu\ge 3\nu_2>\frac32n,
$$
i.e., $\mu>\frac{n}2$, a contradiction.

\medskip

\noindent{\it The case D).} As in the case E), we first blow-up $l_0^1$, but
now $B_1^2\cap E_2\in s_2$. Set $\mu'=\mult_{s_2}D^2$, $\mu'\le\mu$.
Then,
$$
{\mc N}_{s_2|U_2}\cong {\mc O}(-1)\oplus {\mc O}(-2),
$$
and let $\phi_3:U_3\to U_2$ be the blow-up of $s_2$, $E_3$ the exceptional divisor
with the minimal section $s_3$ and a fiber $f_3$. We obtain the situation
like in the previous case: $B_1^3\cap T_U^3=\emptyset$, $T_U^3|_{E_3}=s_3$, so
$D^3|_{E_3}$ has the multiplicity not less than $\nu_2$ at the point $B_1^3\cap E_3$.
Since
$$
D^3|_{E_3}\sim \mu' s_3+(2\mu'+\mu+n-2\nu_1)f_3,
$$
and $\mu'\le\mu$, we get the inequality $3\mu+n-2\nu_1\ge\nu_2$, i.e., the same as
in the case E).

\bigskip

\noindent{\bf The cuspidal case.} Let $B_0$ be the cuspidal point of $l_0\in{\mc P}$.
As it is shown in section \ref{sec3}, we can define $U$ locally by
$$
w^2+z^3+zf_4(x,y)+x+f_6(x,y)=0,
$$
where $\deg f_i\le i$, $f_i(0,0)=0$, $f_6$ has no the linear part, 
$B_0=(0,0,0,0)$, $l_0=\{x=y=0\}$, and $T_U=\{x=0\}$. We see that $T_U$ has at $B_0$
one of the following du Val singularities: $A_1$, $A_2$, $D_4$, $E_6$, $E_7$, or 
$E_8$. 

As in the nodal case, let $\phi_1:U_1\to U$ be the blow-up of $B_0$, $E_1$ the
exceptional divisor, and $P=l_0^1\cap E_1$. The strict transforms of divisors from 
$|H-l_0|\setminus\{T_U\}$ cut out a line $e_1$ on $E_1$. Then, $T_U^1$ is non-singular
along $l_0^1$, and $T_U^1|_{E_1}$ is one of the following:
\begin{itemize}
\item a non-singular conic that is tangent to $e_1$ at $P$, if $B_0$ is $A_1$-point
       on $T_U$;
\item a couple of lines that are different from $e_1$ and have their common point at 
       $P$, in the case $A_2$;
\item a double line that is different from $e_1$ and pass through $P$, in the other 
      cases.
\end{itemize}

We have the following possibilities:
\begin{itemize}
\item[\it A).] $B_1\ne e_1$ and $P\not\in B_1$;
\item[\it B).] $B_1=e_1$;
\item[\it C).] $P\in B_1$ and $B_1\subset\Supp T_U^1|_{E_1}$;
\item[\it D).] $P\in B_1\ne e_1$, $B_1\not\subset\Supp T_U^1|_{E_1}$, and 
                $T_U^1|_{E_1}$ consists of lines;
\item[\it E).] $P\in B_1\ne e_1$ and $T_U^1|_{E_1}$ is a non-singular conic.
\end{itemize}
The case A) was excluded in \cite{Grin1}. The argumentation of the cases 
B)--E) of the nodal situation works in the caspidal situation as well, we 
only give some remarks. One keeps the notation of the corresponding cases.

First, there are no changes in the cases C), D), and E). Then, in the case B),  
let $Q=l_0^2\cap E_2$; $Q$ is not contained in the minimal section. Elements
from ${\mc L}^2|_{E_2}$ (except for $T_U^2$) are all tangent to each other
at $Q$. Suppose that $B_2$ is also tangent to them at $Q$ (if not, the argumentation
of the case $k<N$ works). Then, for a general $S\in{\mc L}$, $S^3$ cuts out
a double fiber of $E_3$, along which it has a double point. We can always
assume that $B_3$ does not pass through the singular point of $S^3$.
Since $B_3$ is also tangent to $S^3$, we will have the analoguos
situation on $U_4$, and so on. Anyway, we may always assume $S^i$ to be non-singular
at the point $\tilde B_i=S^i\cap B_i$. There are no other changes in this case.

\bigskip

So, taking into account the results of \cite{Grin1}, 
we proved the following proposition:
\begin{proposition}
\label{U_prop}
Suppose that a linear system ${\mc D}$ on $U$ has no maximal curves. Then 
${\mc D}$ has no maximal singularities over points.
\end{proposition}

\section{Excluding super-maximal singularities over non-singular points.}
\label{sec5}

In this section we assume the reader to be more or less familiar with
the "super-maximal" modification of the method of maximal singularities
(see \cite{Pukh2}). All denotations are close to those in \cite{Grin2}.

Let $V$ be the blow-up of a curve from ${\mc P}$ on $U$, $\rho:V\to\POn$ the 
corresponding fibered structure, $s_0$ and $f$ are the classes of effective
generators of 1-dimensional cycles on $V$, and $\{G_V\}=|-K_V-2F|$ (see section
\ref{sec3}).   

Consider a linear system ${\mc D}\subset|n(-K_V)+mF|$, $m\ge 0$, without fixed
components, and general elements $D_1, D_2\in{\mc D}$. Then
$$
D_1\circ D_2\sim n^2s_0+(3n^2+2mn)f,
$$
and we set
$$
D_1\circ D_2=Z^h+\sum_{t\in\POn}Z_t^v,
$$
where $Z^h$ and $Z_t^v$ are horizontal and vertical (over $t\in\POn$) 
effective 1-cycles.

We recall that ${\mc D}$ has no maximal singularities over curves (\cite{Pukh2},
\cite{Grin2}). Assume ${\mc D}$ has a super-maximal singularity. Then we have
\begin{equation}
\label{sm_ineq}
m<\sum_{t\in\POn} 
\max_{{\mf v}\in{\mc M}\atop \Center\: {\mf v}\in F_t} 
\frac{e({\mf v})}{{\mf v}(F_t)},
\end{equation}
where $e({\mf v})={\mf v}({\mc D})-\delta_{\mf v}$,
$\delta_{\mf v}$ is the canonical multiplicity with respect to 
${\mf v}$, and ${\mc M}=\{{\mf v}: e({\mf v})>0\}$ 
is the (finite) set of maximal valuations of ${\mc D}$.

Let $A$ be the number of elements in ${\mc M}$, and define 
$\deg Z^h=Z^h\circ F$, $\deg Z_t^v=2Z_t^v\circ (-K_V)$. Then
$$
\sum_{t\in\POn}\deg Z_t^v\le 6n^2+4mn<
\sum_{t\in\POn}\left(\frac{6n^2}{A}+
4n\max_{{\mf v}\in{\mc M}\atop \Center\: {\mf v}\in F_t} 
\frac{e({\mf v})}{{\mf v}(F_t)}\right).
$$
Thus, there exist $t\in\POn$ and ${\mf v}\in{\mc M}$ such that
$$
\deg Z_t^v<\frac{6n^2}{A}+4n\frac{e({\mf v})}{{\mf v}(F_t)}.
$$ 
Consider the graph $\Gamma=\Gamma({\mf v})$ (see section \ref{sec2}). Suppose 
it consists of $N$ vertices, the vertices $\{1,\ldots, L\}$ correspond to points 
(i.e., $B_0,\ldots,B_{L-1}$ are points), and the others correspond to curves.
We denote by $r_i$ the number of all paths in $\Gamma$ from $N$ to $i$, and define
$\Sigma_0=\sum_{i=1}^Lr_i$, $\Sigma_1=\sum_{i=L+1}^Nr_i$.
Let
$$
L'=\max\{1\le i\le L: B_{i-1}\in F_t^{i-1}\}
$$
and $\Sigma'_0=\sum_{i=1}^{L'}r_i$; clearly, $\Sigma'_0\le\Sigma_0$. Then 
$\Sigma'_0\le{\mf v}(F_t)$, and we have
$$
\Sigma'_0 \deg Z_t^v<\frac{6n^2}{A}+4ne({\mf v}).
$$
Put $m_i^h=\mult_{B_{i-1}}\left(Z^h\right)^{i-1}$ and
$m_i^v=\mult_{B_{i-1}}\left(Z^v_t\right)^{i-1}$. Then we have
the quadratic inequality
$$
\sum_{i=1}^Lr_im_i^h+\sum_{i=1}^{L'}r_im_i^v\ge
\frac{\left(2n\Sigma_0+n\Sigma)1+e({\mf v})\right)^2}{\Sigma_0+\Sigma_1}
\ge 4n^2\Sigma_0+4ne({\mf v}).
$$
Note that $m_i^h\le\deg Z^h\le n^2$. Taking into account (\ref{sm_ineq}), we obtain
$$
n^2\Sigma_0+\frac{6n^2}{A}\Sigma'_0>4n^2\Sigma_0,
$$
hence
$$
\frac2A>\frac{\Sigma_0}{\Sigma'_0}\ge 1,
$$
which yields $A=1$ and $\Sigma'_0>\frac12\Sigma_0$.

Thus, it is shown that ${\mc D}$ only consists of valuations centered over points
in $F_t$. This allows us to re-write the super-maximal condition for 
${\mf v}$ in the form
$$
e({\mf v})>m{\mf v}(F_t).
$$ 
Let $B_0=\Center_V{\mf v}$. We change the denotation $F_t$ for the 
"central" fiber (i.e., that contains the super-maximal singularity), denoting
it by $F_0$. The anticanonical linear system on $F_0$ contains a unique element,
say, $l_0$, that passes through $B_0$. 

We summarize the known results as follows (we assume $B_0$ to be a non-singular point
of $V$):
\begin{itemize}
\item ${\mc D}$ may only have maximal singularities over points in $F_0$, and 
       at least one of the maximal singularities is super-maximal;
\item $B_0$ is necessarily contained in a section of the class $s_0$, thus
      $B_0\in G_V$ and $l_0=G_V\cap F_0$ (\cite{Grin2}, \S 2.2);
\item $B_0$ is a non-singular point of $F_0$, $B_0\ne \Bas|-K_{F_0}|$, and
        $l_0$ has a double point (node or cusp) at $B_0$ (\cite{Grin2},
         lemma 2.8);
\item we have the following inequality for ${\mf v}$:
        $\sum_{i=1}^Nr_i\nu_i>2n\Sigma_0+n\Sigma_1+m\Sigma'_0$, where
        $\nu_i=\mult_{B_{i-1}}{\mc D}^{i-1}$;
\item always $\Sigma'_0>\frac12\Sigma_0$.
\end{itemize}
Note that if $V\to U$ is the blow-up of a curve $C\in{\mc P}$, then 
$l_0\cong C$ and $G_V=\POn\times C$. Whence $C$, and thus $l_0$
and $V$ itself, must be singular.

We can say a bit more about the situation. Since $B_0=\Center_{V}{\mf v}$,
${\mc M}$ has to contain a discrete valuation that defines a maximal singularity
over $B_0$ and can be realized by a divisorial extraction (i.e., a wheighted 
blow-up, see section \ref{sec2}). We denote such a discrete valuation by 
${\mf v}_{div}$.

We have three possible cases:
\begin{itemize}
\item[\bf A.] ${\mf v}_{div}$ is super-maximal (for example,
   ${\mf v}_{div}={\mf v}$) and $B_0$ is a 
     non-singular point of $V$;
\item[\bf B.] ${\mf v}_{div}$ is not super-maximal (thus 
      ${\mf v}_{div}\ne{\mf v}$), $B_0$ is non-singular;
\item[\bf C.] $B_0$ is the singular point of $V$.
\end{itemize}
The latter case is studied in the next section. Here we show that the first two cases
are impossible. One can always assume $e({\mf v}_{div})\ge e({\mf v})$, 
so in the case {\bf B} ${\mf v}_{div}(F_0)>{\mf v}(F_0)$, and we have the 
weak super-maximal condition for ${\mf v}_{div}$:
$$
e({\mf v}_{div})>m{\mf v}(F_0).
$$

\noindent{\bf The case A: ${\mf v}_{div}$ is super-maximal.} We may assume
that ${\mf v}_{div}={\mf v}$.

Suppose first that $B_0$ is a nodal point of $l_0$. As it is mentioned in section 
\ref{sec2}, ${\mf v}$ can be realized as the weighted blow-up with weights
$(1,L,N)$, $L<N$, or as a chain of blow-ups with centers $B_0,\ldots,B_{N-1}$, where
$B_0,\ldots,B_{L-1}$ are points and $B_L,\ldots,B_{N-1}$ are curves, with the
conditions $B_i\not\subset E_{i-1}^i$ for all $i$ and $B_i\cap E_{i-1}^i=\emptyset$ 
for $i>L$. One uses the latter realization. We set
$$
L'=\max\{1\le i\le L: B_{i-1}\in F_0^{i-1}\}
$$
and
$$
k=\max\{1\le i\le L: B_{i-1}\in l_0^{i-1}\}.
$$
Clearly, $L'\le{\mf v}(F_0)$ and $k\le L'$. Since $r_i=1$ for all $i$, we 
have the super-maximality condition in the form
$$
\sum_{i=1}^N\nu_i >n(N+L)+mL'.
$$
Choose general elements $D_1,D_2\in{\mc D}$, and let $D_1\circ D_2=Z^h+Z^v$, where
$Z^h$ and $Z^v$ are horizontal and vertical 1-cycles. $Z^v$ contains the 1-cycle
$Z_0^v$ that lies in $F_0$, and let
$$
Z_0^v=\alpha l_0+\tilde Z_0^v,
$$ 
where $l_0\not\subset\Supp\tilde Z_0^v$. Define the multiplicities
$m_i^h=\mult_{B_{i-1}}(Z^h)^{i-1}$ 
and $\tilde m_i^v=\mult_{B_{i-1}}(\tilde Z_0^v)^{i-1}$ for $1\le i\le L$.
We obtain the quadratic inequality in the form
$$
\sum_{i=1}^Lm_i^h+\alpha(k+1)+\sum_{i=1}^{L'}\tilde m_i^v>
\frac{\left(n(L+N)+mL'\right)^2}{N}\ge 4n^2L+4mnL'.
$$
Since $m_i^h\le\deg Z^h\le n^2$ and 
$\alpha+\tilde m_i^v\le\deg Z^v\le 3n^2+2mn$, we deduice
$$
n^2L+\alpha+3n^2L'+2mnL'>4n^2L+4mnL',
$$
whence
$$
\alpha>3n^2(L-L')+2mnL'.
$$
But $\alpha\le 3n^2+2mn$, hence $L=L'$. Moreover, the same arguments show that
$k=L'$. Thus, in what follows, we suppose that $k=L'=L$. In particular, 
$B_i\cap F_0^i\ne\emptyset$ for all $i$.

Now we show that $B_L\not\subset F_0^L$. Assume the converse. Set
$e_i=E_i\cap F_0^i$, and choose a point $\tilde B_L\in e_L$ outside
of $e_{L-1}^L\cup l_0^L$. Clearly, 
$\tilde\nu_{L+1}\eqdef\mult_{\tilde B_L}(D^L|_{F_0^L})\ge\nu_{L+1}$.
For a general $D\in{\mc D}$ we set
$$
D|_{F_0}=\mu l_0+C\in|n(-K_{F_0})|,
$$
where $\mu\le n$ and $l_0\not\subset\Supp C$. Define $\alpha_i=\mult_{B_i}C^i$
for $i<L$, $\tilde\alpha_L=\mult_{\tilde B_L}C^L$, 
$m_i=\ord_{e_i}(D^i|_{F_0^i})$, and 
$\tilde\nu_i=\mult_{B_{i-1}}(D^{i-1}|_{F_0^{i-1}})$. Then we have
$$
\begin{array}{l}
\tilde\nu_1=2\mu+\alpha_0, \\
\tilde\nu_2=\mu+\alpha_1+m_1, \\
\cdots \\
\tilde\nu_L=\mu+\alpha_{L-1}+m_{L-1}, \\
\tilde\nu_{L+1}=\tilde\alpha_L+m_L.
\end{array}
$$
Since 
$$
\alpha_0+\ldots+\alpha_{L-1}+\tilde\alpha_L\le 2\alpha_0+\alpha_1+\ldots
+\alpha_{L-1}\le (D|_S-\mu l_0)\circ l_0=n-\mu
$$
and $\tilde\nu_i=\nu_i+m_i$ for $1\le i\le L$, we obtain
$$
\nu_1+\ldots+\nu_{L+1}\le (L+1)\mu+n-\mu\le (L+1)n,
$$
which contradicts to the condition $\nu_i>n$ for all $i$.

Thus, we assume that $B_L\not\subset F_0^L$. We define the points 
$\tilde B_i=B_i\cap F_0^i$, $0\le i\le N-1$. Suppose that 
$\tilde\nu_i=\mult_{\tilde B_{i-1}}(D^{i-1}|_{F_0^{i-1}})$, $1\le i\le N$,
the numbers $m_i$, $\alpha_i$, $\mu$ are defined as before. Set
$$
k'=\max\{1\le i\le N: \tilde B_{i-1}\in l_0^{i-1}\}.
$$
Clearly, $k'\ge k=L'=L$. Assume that $k'<N$. Then we have
$$
\begin{array}{l}
\tilde\nu_1=2\mu+\alpha_0, \\
\tilde\nu_2=\mu+\alpha_1+m_1, \\
\cdots \\
\tilde\nu_{k'}=\mu+\alpha_{k'-1}+m_{k'-1}, \\
\tilde\nu_{k'+1}=\tilde\alpha_{k'}+m_{k'},
\end{array}
$$
and since 
$$\alpha_0+\ldots+\alpha_{k'}\le 2\alpha_0+\alpha_1+\ldots+\alpha_{k'-1}\le n-\mu,
$$
we get a contradiction:
$$
\nu_1+\ldots+\nu_{k'+1}\le (k'+1)\mu+n-\mu\le (k'+1)n.
$$
Thus, $k'=N$. We repeat the arguments to obtain
$$
\nu_1+\ldots+\nu_N\le (N+1)\mu+n-\mu\le (N+1)n,
$$
but this is impossible because of the super-maximality condition
$$
\nu_1+\ldots+\nu_N>n(L+N)+mL'\ge (N+1)n.
$$
The important remark is that we have used the N\"other-Fano inequality,
not the super-maximality condition.

So, $B_0$ has to be a cusp of $l_0$. It is easy to observe that nothing
changes in the cuspidal case, the above argumentation does work. In fact,
this case is even simpler since $B_2$ avoids $l_0^2$ (note that 
$l_0^2\cap E_1^2\ne\emptyset$).

\medskip

\noindent{\bf The case B: ${\mf v}_{div}$ is not super-maximal.}
Since ${\mf v}_{div}$ is not super-maximal, one has 
$e({\mf v}_{div})\le m{\mf v}_{div}(F_0)$ . We may always assume that 
$e({\mf v}_{div})>e({\mf v})$.
Since ${\mf v}$ is super-maximal, that is, $e({\mf v})>m{\mf v}(F_0)$, we
obtain two important inequalities:
$$
\begin{array}{l}
e({\mf v}_{div})>m{\mf v}(F_0), \\
{\mf v}_{div}(F_0)>{\mf v}(F_0).
\end{array}
$$
Suppose that ${\mf v}_{div}$ is realized by the weighted blow-up with weights
$(1,L,N)$. Taking into account the scheme of resolution of weighted blow-up,
we have the centers $B_0,\ldots,B_{N-1}$, where $B_0,\ldots,B_{L-1}$ are points
and $B_L,\ldots,B_{N-1}$ are curves. Set
$$
k=\max\{1\le i\le L: B_{i-1}\in l_0^{i-1}\}.
$$
We repeat the arguments of the case {\bf A} to argue that 
$B_L\not\subset F_0^L$. Moreover, due to the remark at the end of the case
{\bf A}, that is, that we have not needed the super-maximality condition,
we also argue that $k<L$.
We show that $B_k\not\in F_0^k$ for $k<L$. Assume the converse, that is,
$B_k\in F_0^k\setminus l_0^k$. Restricting ${\mc D}$ to $F_0$ and 
using the denotations of $\tilde\nu_i$, $\mu$, $\alpha_i$, and so on,
in the same way as above, we obtain
$$
\begin{array}{l}
\tilde\nu_1=2\mu+\alpha_0, \\
\tilde\nu_2=\mu+\alpha_1+m_1, \\
\cdots \\
\tilde\nu_k=\mu+\alpha_{k-1}+m_{k-1}, \\
\tilde\nu_{k+1}=\alpha_k+m_k,
\end{array}
$$
whence
$$
\nu_1+\ldots+\nu_{k+1}\le (k+1)\mu+\alpha_0+\ldots+\alpha_k\le
(k+1)\mu+n-\mu\le (k+1)n,
$$
which is impossible. Thus we have
$$
{\mf v}_{div}(F_0)=k<L,
$$
and since ${\mf v}_{div}(F_0)>{\mf v}(F_0)\ge 1$, then $k\ge 2$.

Let $Z_0^v=\alpha l_0+\tilde Z_0^v$, and let $m_i^h$ and $\tilde m_i^v$
be defined as above. One uses the weak super-maximality
condition
$$
\sum_{i=1}^N\nu_i>n(L+N)+m{\mf v}(F_0)
$$
to obtain the quadratic inequality in the form
$$
\sum_{i=1}^Lm_i^h+(k+1)\alpha+\sum_{i=1}^k\tilde m_i^v>
\frac{\left(n(L+N)+m{\mf v}(F_0)\right)^2}{N}.
$$
Remark that 
$$
2\tilde m_1^v+\tilde m_2^v+\ldots\tilde m_k^v\le 
\tilde Z_0^v\circ l_0\le 3n^2+2mn-\alpha,
$$
and the more so
$$
\tilde m_1^v+\tilde m_2^v+\ldots\tilde m_k^v\le 3n^2+2mn-\alpha.
$$
Moreover,
$$
m_1^h+\ldots+m_k^h\le Z^h\circ F_0\le n^2.
$$
Thus we get
$$
\frac{L}{k}n^2+k\alpha+3n^2+2mn>
4n^2L+4mn{\mf v}(F_0)+n^2\frac{\left((N-L)-\frac{m}{n}{\mf v}(F_0)\right)^2}{N},
$$
hence
$$
k\alpha>n^2\left(4L-\frac{L}{k}-3\right)+2mn\left(2{\mf v}(F_0)-1\right)+
   n^2\frac{\left((N-L)-\frac{m}{n}{\mf v}(F_0)\right)^2}{N}.
$$
We observe that $\nu_1+\ldots+\nu_k\le (k+1)\mu+\alpha_0+\ldots+\alpha_{k-1}\le
(k+1)n$, and then
$$
\theta_k\eqdef\frac{\nu_1+\ldots+\nu_k}{kn}\le\frac{k+1}{k}.
$$
Clearly, 
$$
n\theta_k\ge\frac1N\sum_{i=1}^N\nu_i,
$$
and we deduice from the weak super-maximal condition that
$$
N-L>\frac{2-\theta_k}{\theta_k-1}L+
\frac1{\theta_k-1}\cdot\frac{m}{n}{\mf v}(F_0).
$$
Combining this and the estimation for $\theta_k$ with the quadratic
inequality, one gets
$$
\alpha>n^2\left(\frac{(k+2)L}{k}-\frac3k\right)+
  mn\left(\frac{(k+1)^2}{k^2}{\mf v}(F_0)-\frac2k\right).
$$
Recall that $L\ge k+1$. Thus we have
$$
\alpha>n^2\left(\frac{k^2+3k+1}{k}\right)+
mn\left({\mf v}(F_0)+\frac{2{\mf v}(F_0)}{k}+
\frac{{\mf v}(F_0)}{k^2}-\frac2k\right),
$$
and the condition $\alpha\le 3n^2+2mn$ yields
$$
{\mf v}(F_0)=1.
$$
Note that the argumentation is the same in the nodal and in the cuspidal cases.

Now we turn back to the super-maximal valuation ${\mf v}$. All denotations 
($\nu_i$, $B_i$, $E_i$, and so on) are related to ${\mf v}$. Since 
$\Sigma'_0\le{\mf v}(F_0)=1$ and $\Sigma_0<2\Sigma'_0$, we have
$\Sigma_0=\Sigma'_0=1$. Thus $L=1$ and $B_1,\ldots,B_{N-1}$ are curves.
It is known that $B_{L+1}\cap E_L^{L+1}=\emptyset$ and 
$B_i\not\subset E_{i-1}^i$ for all $i>L$. In our case $L=1$, thus 
$B_2\cap E_1^2=\emptyset$, $E_2\cong\FA_2$. Denote by $s_i$ and $f_i$ the 
minimal sections and fibers of the corresponding ruled surfaces $E_i$. We see
that $B_2\sim s_2+2f_2$ and $B_i\ne s_i$ for all $i>1$. Set
$$
l=\max\{3\le i\le N: B_{i-1}\cap E_{i-2}^{i-1}=\emptyset\}.
$$
Assume that $l<N$. Then $B_l\sim s_l+bf_l$, where $b>l$. Indeed, from lemma
\ref{Ei_lem} for the case $L=1$ it follows that $E_i\cong\FA_i$ for $i\le l$.
Take into account that $E_{l-1}^l|_{E_l}=s_l$, thus $B_l\cap s_l\ne\emptyset$
and $b>l$. It is easy to compute that
$$
D^l|_{E_l}\sim \nu_ls_l+(\nu_1+\ldots+\nu_l)f_l.
$$
Since $\nu_{l+1}B_l\subset D^l|_{E_l}$, we get
$$
(l+1)\nu_{l+1}\le b\nu_{l+1}\le\nu_1+\ldots+\nu_l.
$$
Set
$$
k=\max\{1\le i\le N-1: B_{i-1}\cap l_0^{i-1}\ne\emptyset\}.
$$
Assume that $k<l$. Restricting ${\mc D}$ to $F_0$, we have
$$
\begin{array}{l}
\tilde\nu_1=2\mu+\alpha_0, \\
\tilde\nu_2=\mu+\alpha_1+m_1, \\
\cdots \\
\tilde\nu_{k+1}=\alpha_k+m_k.
\end{array}
$$
Since $2\alpha_0+\alpha_1+\ldots+\alpha_{k-1}\le n-\mu$, we get
$$
\nu_1+\ldots+\nu_{k+1}\le (k+1)\mu+n-\mu\le (k+1)n,
$$
which is impossible because of the inequalities $\nu_i>n$ for all $i$.
Thus $k\ge l$ and
$$
\tilde\nu_l=\mu+\alpha_{l-1}+m_{l-1}.
$$
We obtain
$$
\nu_1+\ldots+\nu_l\le (l+1)\mu+n-\mu\le (l+1)n.
$$
Combining this estimation with $(l+1)\nu_{l+1}\le\nu_1+\ldots+\nu_l$, we get
a contradiction:
$$
\nu_{l+1}\le n.
$$
Thus, we obtain $l=N$. In other words, the valuation ${\mf v}$ corresponds to
a weighted blow-up. In the case {\bf A} it is shown that this is impossible.
We proved the following result:
\begin{proposition}
\label{V_prop}
The linear system ${\mc D}\subset|n(-K_V)+mF|$, $m\ge 0$, has no super-maximal
singularities over the non-singular part of $V$.
\end{proposition}

\section{Excluding super-maximal singularities over the singular point.}
\label{sec6}

It remains to exclude super-maximal singularities over the point $B_0$ which
is the singular point of $V/\POn$. Let ${\mf v}$ be a maximal valuation with
the center at $B_0$ which is a super-maximal at the same time, that is, the
inequality
$$
e({\mf v})>m{\mf v}(F_0)
$$
holds, where $F_0$ is the fiber that contains $B_0$. As before, there exists
a maximal valuation ${\mf v}_{div}$ which is a divisorial extraction in the
Mori category, $\Center_V{\mf v}_{div}=B_0$. We may always assume that 
$e({\mf v}_{div})\ge e({\mf v})$. Hence ${\mf v}_{div}$ satisfies the weak
super-maximal condition $e({\mf v}_{div})>m{\mf v}(F_0)$. In particular,
one has always
$$
e({\mf v}_{div})>m.
$$
As it follows from Kawakita's theorem \ref{Kaw_th}, only the following cases are 
possible:
\begin{itemize}
\item $B_0$ is the simplest double point, that is, it is locally defined by
   $xy+zw=0$, and ${\mf v}_{div}$ is realized by the blow-up of $B_0$;
\item $B_0$ is defined by $xy+z^3+w^2=0$ and ${\mf v}_{div}$ is realized by
  the blow-up
  of $B_0$;
\item $B_0$ is defined by $xy+z^3+w^2=0$ and ${\mf v}_{div}$ is realized by 
the weighted blow-up with weights $(1,5,2,3)$. 
\end{itemize}
Now we argue that the first two cases are impossible. Let $\gamma_1:U_1\to U$
be the blow up of $B_0$, $E_1$ the exceptional divisor. In the first two cases
$E_1$ realizes ${\mf v}_{div}$, and either $E_1\cong\POn\times\POn$ if $B_0$ is
the simplest double point, or $E_1$ is isomorphic to a quadric cone in $\PT$ (in
the second case). Remark that $U_1$ is non-singular. For a general 
$D\in{\mc D}\subset|n(-K_V+mF)|$ we set
$$
D^1=\gamma_1^*(D)-\nu_1E_1.
$$
The N\"other-Fano inequality has the form $\nu_1>n$. Choose a general curve
$C\in|2(-K_{F_0})-B_0|$. Then $C$ is singular at $B_0$ and $C^1\circ E_1=2$.
We get a contradiction immediately:
$0\le D^1\circ C^1=2n-2\nu_1<0$. 

In what follows we assume that $B_0$ is locally defined by $xy+z^3+w^2=0$, 
${\mf v}_{div}$ corresponds to the third case and satisfies 
$e({\mf v}_{div})>m$. 

\begin{figure}[htbp]
\begin{center}
\epsfxsize 10cm
\hepsffile{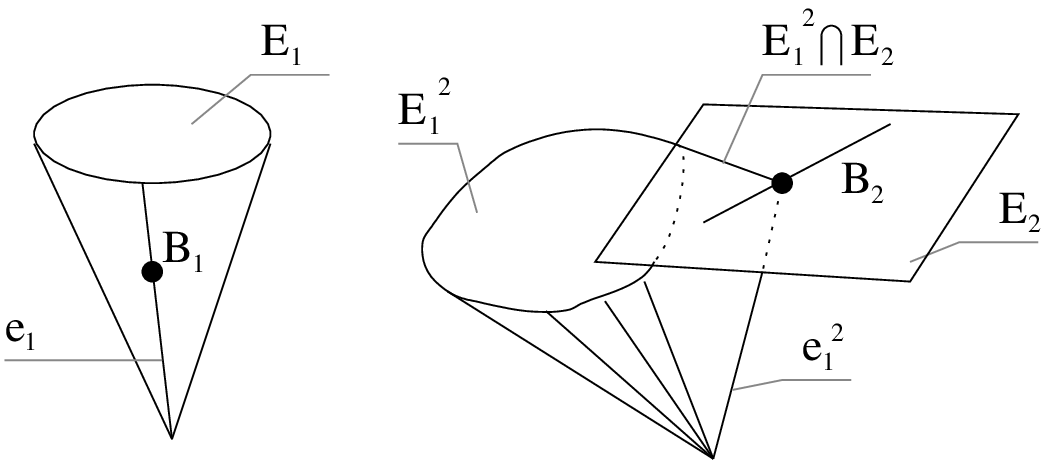}
\caption{}\label{ver2.fig}
\end{center}
\end{figure}

The valuation ${\mf v}_{div}$ can be geometrically realized as follows 
(see \cite{Kaw2}, lemma 5.2). Let $\psi_1:V_1\to V$ be the blow-up of the point
$B_0$. The exceptional divisor $E_1$ is isomorphic to a quadric cone in $\PT$.
The center of ${\mf v}_{div}$ on $V_1$ is a point $B_1\in E_1$ lying outside
of the vertex of the cone. By $e_1\subset E_1$ we denote the generator of the 
cone that passes through $B_1$ (see figure \ref{ver2.fig}). Let 
$\psi_2:V_2\to V_1$ be the blow-up of $B_1$, $E_2\cong\PTw$ the exceptional
divisor. Then $B_2=\Center_{V_2}{\mf v}_{div}$ is a line in $E_2$ that passes
through the point $e_1^2\cap E_2$ and does not contained in $E_1^2$. The 
exceptional divisor $E_3$ of the blow-up $\psi_3:V_3\to V_2$ of $B_2$ 
realizes ${\mf v}_{div}$.

We set $\nu_i=\mult_{B_{i-1}}{\mc D}^{i-1}$. Observe that $2\nu_1\ge\nu_2\ge\nu_3$.
We assume that $\nu_1\le n$, otherwise we get the situation of the first two cases.
The weak super-maximal condition for ${\mf v}_{div}$ has the form
$$
\nu_1+\nu_2+\nu_3>4n+m.
$$

\begin{figure}[htbp]
\begin{center}
\epsfxsize 12cm
\hepsffile{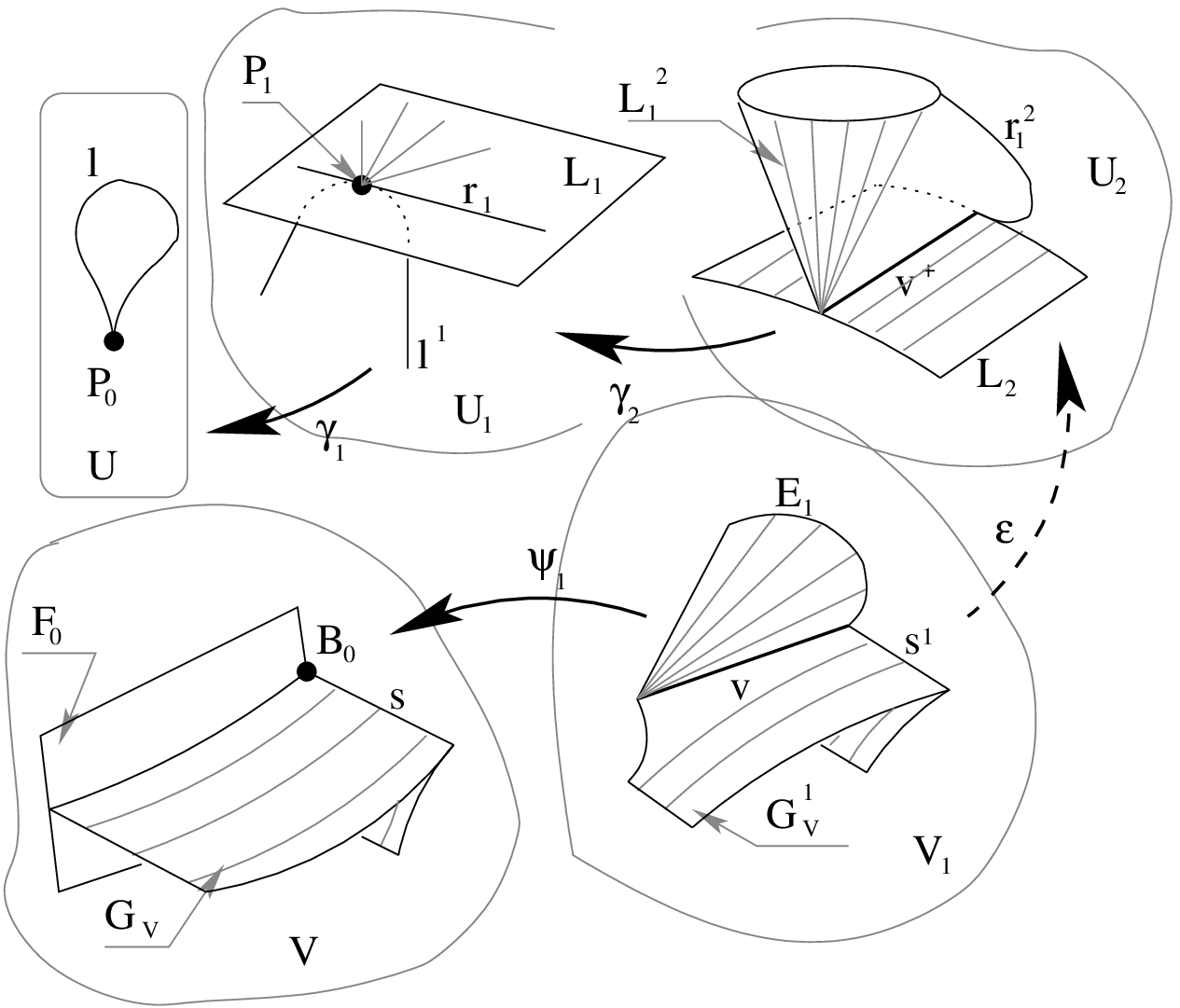}
\caption{}\label{ver3.fig}
\end{center}
\end{figure}

We assume that $V$ is obtained by the blow-up $\phi:V\to U$ of a curve 
$l\in{\mc P}$ on $U$, and $l$ has a cusp at a point $P_0$. The birational
morphism $\phi$ can be also realized as follows (see figure \ref{ver3.fig}).
Let $\gamma_1:U_1\to U$ be the blow-up of $P_0$, $L_1\cong\PTw$ the exceptional
divisor. The strict transform $l^1$ of the curve $l$ on $V_1$ is tangent
to $L_1$ at a point $P_1$. The tangent direction to $l^1$ at $P_1$ defines
a line $r_1\subset L_1$. Let $\gamma_2:U_2\to U_1$ be the blow-up of $l^1$,
$L_2$ the exceptional divisor. The strict transform $L_1^2$ of $L_1$ on $V_2$
becomes isomorphic to a quadric cone that is blown up at some point outside of 
the vertex of the cone, and $r_1^2$ is the $(-1)$-curve of $L_1^2$. Note that
$L_1^2$ is tangent to $L_2$ along the fiber $v^+$.

Let us turn back to $V$ for instance. We recall that the exceptional divisor
$G_V$ of the morphism $\phi$ is isomorphic to $l\times\POn$ (see section 
\ref{sec3}). Hence $G_V$ is singular along a section $s$ of $V\to\POn$. This
section has the class $s_0$. Now let $\psi_1:V_1\to V$ be the blow-up of $B_0$, 
$E_1$ the exceptional divisor. As it is mentioned above, $E_1$ is isomorphic 
to a quadric cone in $\PT$. The divisor $G_V^1$ is tangent to $E_1$ along
the generator $v$. Remark that $s^1$ intersects $E_1$ outside of the vertex of the
cone. 

Finally, there exists a flop $\eps:V_1\dasharrow U_2$ centered at $s^1$. We 
observe that the flop $\eps^{-1}$ is centered at $r_1^2$. The birational map
$\eps$ maps $E_1$ onto $L_1^2$, $v$ onto $v^+$, and $G_V^1$ onto $L_2$.

It is worth to give an explicite equation for $V$ in a neighborhood of the point $B_0$.
As in section \ref{sec3}, we assume that $U$ is defined by the equation
$$
w^2+z^3+zf_4(x,y)+x+f_6(x,y)=0
$$
in a neighborhood of the point $P_0=(0,0,0,0)$, the polynomials $f_i$ satisfy the
same conditions as before. The blow-up $\phi$ is given by $x=ty$, so we have for $V$
$$
w^2+z^3+zf_4(ty,y)+ty+f_6(ty,y)=0.
$$
We observe that $B_0=(0,0,0,0)$, $G_V=\{y=0\}$, the fiber $F_0=\{t=0\}$, and
$s=\{w=z=y=0\}$. 

We assume that ${\mc D}_V={\mc D}\subset|n(-K_V+mF)|$. Denote by ${\mc D}_U$
the strict transform of ${\mc D}_V$ on $U$, by ${\mc D}_{U_1}$ the strict
transform on $U_1$, and so on. Suppose that ${\mc D}_U\subset|aH|$, and let
$\mu=\mult_l{\mc D}_U$ and $\nu^*=\mult_{P_0}{\mc D}_U$. It is easy to compute
that
$$
\begin{array}{l}
a=2n+m, \\
\mu=n+m, \\
\nu^*=\nu_1+\mu=\nu_1+n+m.
\end{array}
$$
There are the following possibilities:
\begin{itemize}
\item[\bf A.] $e_1=v$ and $B_1=s^1\cap E_1$ (see figures \ref{ver2.fig} and
       \ref{ver3.fig});
\item[\bf B.] $e_1=v$ and $B_1\not\in s^1$;
\item[\bf C.] $B_1\in E^1\setminus G_V^1$ (that is, $e_1\ne v$).
\end{itemize}
We deal with them separately.

\medskip

\noindent{\bf The case A: $e_1=v$ and $B_1=s^1\cap E_1$.} The divisor $G_V$
is cut out by the hypersurface $\{y=0\}$. In the local coordinates $[w,z,t]$
the equation $\{w^2+z^3=0\}$ defines $G_V$. According to the construction
of $V$ via double covering, we see that $G_V$ is a double cover of a quadric
$G\cong\POn\times\POn$ with the ramification divisor $3p+p'$, where $p$ and $p'$
have the class $(1,0)$ on $G$ and $p\ne p'$. The quadric $G$ has the local 
coordinates $[z,t]$. We assume that $p$ is defined by the equation $p=\{z=0\}$. 
Remark that $s$ lies over $p$. Consider a curve $\tilde C\subset G$ of the class
$(1,2)$ defined by the equation $\{z=at^2\}$ for some $a\ne 0$. Let $C$ be the
pre-image of $\tilde C$ on $G_V$. Thus, $C$ is locally defined by the equations
$\{y=w^2+a^3t^6=0\}$. Using the equation for $V$, it is easy to check that
$C^1$ has a double point at the point $s^1\cap E_1$, and $C^2$ also has a double
point at $s^2\cap E^2$ (see figure \ref{ver2.fig} and note that $e_1=v$ and 
$B_1=s^1\cap E_1$ in the considered case). Taking into account that $C$ has the
class $2s+2f$ on $V$, for a general $D\in{\mc D}_V$ we obtain
$$
D^2\circ C^2=4n+2m-2\nu_1-2\nu_2\ge 0,
$$
hence $\nu_1+\nu_2\le 2n+m$. Since $\nu_3\le 2n$, we get a contradiction with
the super-maximal condition:
$$
4n+m\ge \nu_1+\nu_2+\nu_3>4n+m.
$$

\noindent{\bf The case B: $e_1=v$ and $B_1\not\in s^1$.} For the curves $v$
and $v^+$, we define
$$
\eta=\mult_v{\mc D}_{V_1}=\mult_{v^+}{\mc D}_{U_2}.
$$
Since $B_1\not\in s^1$, the flop $\eps:V_1\dasharrow U_2$ is a local
isomorphism in a neighborhood of $B_1$, and we denote the image of $B_1$
on $U_2$ by $B_1$ to avoid complicating the notation. Thus, the picture
of the blow-up of $B_1$ given by figure \ref{ver2.fig} can be transfered
to $U_2$.

It is easy to compute that ${\mc N}_{v^+|U_2}\simeq{\mc O}\oplus{\mc O}(-1)$.
We blow up the point $B_1$ and then the strict transform of $v^+$, and denote
by $\phi':U'\to U_2$ the composition of these blow-ups. Let $E'$ be the 
exceptional divisor that corresponds to the blow-up of the strict transform of
$v^+$. We see that $E'\cong\FA_1$, and denote by $s'$ and $f'$ the minimal
section and a fiber of $E'$. We have
$$
{\mc D}_{U'}|_{E'}\subset|\eta s'+(2\eta+\mu-\nu_2)f'|.
$$
Moreover, the strict transform of $L_2$ on $U'$ cut out the minimal section
$s'$. Thus, comparing the situation with figure \ref{ver2.fig}, we see
that the strict transform of the curve $B_2$ intersects $E'$ at a point
lying outside of the minimal section. Consequently, ${\mc D}_{U'}|_{E'}$ has
the multiplicity not less than $\nu_3$ at this point. We obtain
$$
2\eta+\mu-\nu_2\ge\nu_3,
$$
or $2\eta+\mu\ge\nu_2+\nu_3$. Since $\nu_1\le n$ and $\nu_1+\nu_2+\nu_3>4n+m$, 
we get $\nu_2+\nu_3>3n+m$. Taking into account that $\mu=n+m$, we obtain
$2\eta>2n$, that is, $\eta>n$. 

Hence the linear system ${\mc D}_{V_1}$ has the multiplicity $\eta>n$ along
one of the generators of the cone $E_1$ (this generator is $v$ in our case),
which yields $E_1\subset\Bas{\mc D}_{V_1}$ because $\nu_1\le n$. 
Thus we get a contradiction.

\medskip

\noindent{\bf The case C: $B_1\not\in E_1\setminus G_V^1$.} In this case 
we have $e_1\ne v$. Consequently, $V_1$ and $U_1$ are isomorphic in a 
neighborhood of the point $B_1$. To avoid complicating the notation,
by $B_1$ and $e_1$ we denote the image of $B_1$ and the strict transform
of $e_1$ on $U_1$ respectively. Note that $e_1$ becomes the line in $L_1$
that passes through $B_1$ and $P_1$, and $e_1\ne r_1$. We set 
$\eta=\mult_{e_1}{\mc D}_{U_1}$ and $\tilde\nu=\mult_{P_1}{\mc D}_{U_1}$.

It is easy to compute that 
${\mc N}_{e_1|U_1}\simeq{\mc O}(-1)\oplus{\mc O}(1)$. On $U_1$, let us first 
blow up the points $B_1$ and $P_1$ and then the strict transform of $e_1$.
We denote by $U'\to U_1$ the composition of these blow-ups, by $E'$ the
exceptional divisor lying over $e_1$, and by $s'$ and $f'$ the minimal section and 
a fiber of $E'$.

We see that $E'\cong\FA_2$, the strict transform of $L_1$ cuts out the minimal
section $s'$ on $E'$. The strict transforms of the curves $l^1$ and $B_2$
intersect $E'$ at two points that lie outside of the minimal section and does
not lie in the same fiber. Thus at these points ${\mc D}_{U'}|_{E'}$ has the 
multiplicities not less than $\mu$ and $\nu_3$ respectively. Taking into
account that
$$
{\mc D}_{U_1}|_{E'}\subset|\eta s'+(3\eta+\nu^*-\tilde\nu-\nu_2)|,
$$ 
where $\nu^*=\nu_1+n+m=\mult_{P_0}{\mc D}_U$, we obtain
$$
3\eta+\nu^*-\tilde\nu-\nu_2\ge\mu+\nu_3.
$$
But $\tilde\nu\ge\mu$, $\nu_2+\nu_3>3n+m$, and $\mu=n+m$, thus we get
$\eta>m$. We get the contradiction exactly as in {\bf the case B}.

\medskip

So, the following proposition is proved:
\begin{proposition}
\label{Vsing_prop}
Linear systems on $V$ have no super-maximal singularities over the singular
point of $V$.
\end{proposition}

\section{Conclusion.}
\label{sec7}
We are ready to finish proving theorem \ref{main_th}. Let $\rho':V'\to S'$
be a Mori fibration, $\chi:U\dasharrow V'$ a birational map. Consider the
complete linear system ${\mc D}'=|n'(-K_{V'})+\rho'^*(L')|$ on $V'$. We assume
that ${\mc D}'$ is very ample and $L'\in\Pic(S')$ is an ample divisor.

Denote by ${\mc D}$ the strict transform of ${\mc D}'$ on $U$. Suppose that
${\mc D}\subset|aH|$. It follows from proposition \ref{U_prop} that either
$K_V+\frac2a{\mc D}$ is canonical, and then $\chi$ is an isomorphism by
\cite{Corti}, theorem 4.2, or ${\mc D}$ has a maximal singularity over
a curve $l\in{\mc P}$. In the latter case, let $\mu=\mult_l{\mc D}$ and
$\phi:V=U_l\to U$ the blow up of $l$. We have 
$$
{\mc D}_V=\phi_*^{-1}{\mc D}\subset|(a-\mu)(-K_V)+(2\mu-a)F|.
$$ 
Observe that $2\mu-a>0$. Propositions \ref{V_prop} and \ref{Vsing_prop} together
with the results of \cite{Grin2}, section 2.2, show that ${\mc D}_V$ has no 
maximal singularities over curves or super-maximal singularities over points.
Thus $a-\mu=n'$, and from \cite{Corti}, theorem 4.2, it follows that the composition
$\chi\circ\phi:V\dasharrow V'$ is birational over the base. Theorem \ref{main_th}
is proved.

Remark that theorem \ref{main_th} completes studying the birational 
classification problem for non-singular Mori fibrations in del Pezzo
surfaces of degree 1. Indeed, theorem 2.6 together with corollary 2.10
in \cite{Grin2}, proposition 2.12 in \cite{Grin2}, and theorem \ref{main_th}
of the present paper exhaust all possible sets of Mori structures
${\mathcal {MS}}$ for varieties of the indicated type. It only remains
to recall that ${\mathcal {MS}}$ is a birational invariant that uniquely defines
the birational class.

\end{document}